\documentclass[12pt]{article}

\normalbaselines
\catcode `\@=11
\@addtoreset{equation}{section}

\def\section{\@startsection {section}{1}{\z@}{-3.5ex plus -1ex minus
     -.2ex}{2.3ex plus .2ex}{\normalsize\bf}}
\def\subsection{\@startsection{subsection}{2}{\z@}{-3.25ex plus -1ex minus
 -.2ex}{1.5ex plus .2ex}{\normalsize\bf}}

\def\thebibliography#1{\section*{References\markboth
  {REFERENCES}{REFERENCES}}\list
  {[\arabic{enumi}]}{\settowidth\labelwidth{[#1]}\leftmargin\labelwidth
  \advance\leftmargin\labelsep
  \usecounter{enumi}}
  \def\newblock{\hskip .11em plus .33em minus -.07em}
  \sloppy
  \sfcode`\.=1000\relax}
 
\catcode `\@=12
\newcommand{\df}{\mbox{$:=$}}
\newtheorem{rem}{Remark}

\newtheorem{thm}{Theorem}

\font\aa=msam10 scaled\magstep1
\font\ab=msbm10 scaled\magstep1
\font\aba=msbm7 
\font\ac=eufm10 scaled\magstep1
\font\ai=eusb10 scaled \magstep 1

\newcommand{\oo}{ 1}
\newcommand{\prend}{\aa \symbol{'003}}

\newcommand{\pre}{\mbox{\prend}}

\newcommand{\gh}{\mbox{\ai H}}

\newcommand{\gz}{\mbox{\ab Z}}
\newcommand{\gr}{\mbox{\ab R}}
\newcommand{\gcp}{\mbox{\ab CP}}
\newcommand{\gc}{\mbox{\ab C}}
\newcommand{\Ka}{K\"ahler}
\newcommand{\Gras}{\mbox{$G_n({\gc}^{m+n})$}}
\newcommand{\cpn}{\mbox{{\aba CP}$^{n}$}}
\newcommand{\got}[1]{{\mbox{\ac{#1}}}}
\newcommand{\mb}[1]{{\mbox{\boldmath{$#1$}}}}

\begin{document}

\vspace*{2.5cm}
\noindent
{ \bf  COHERENT STATES, PHASES AND SYMPLECTIC AREA OF GEODESIC TRIANGLES}
\vspace{1.3cm}\\
\noindent
\hspace*{1in}
%
\begin{minipage}{13cm}
Stefan Berceanu\vspace{0.3cm}\\
Department of Theoretical Physics, Institute of Physics and Nuclear
Engineering,
 P.O. Box-MG-6, Bucharest-Magurele, Romania \\
E-mail: Berceanu@Theor1.ifa.ro
\end{minipage}

\vspace*{0.5cm}

\begin{abstract}
\noindent

On certain manifolds,  the phase   which appears in the scalar product of two
coherent state vectors is twice the symplectic area of the geodesic triangle
determined by the corresponding points on the manifold and the origin of the
system of coordinates. This result is proved for compact Hermitian symmetric
spaces using the generalization via coherent states of the shape invariant  for geodesic triangles
and  re-obtained on the  complex Grassmannian by brute-force calculation.
 
\end{abstract}

\section{\hspace{-4mm}.\hspace{2mm}INTRODUCTION}
The   aim of this paper is to find a geometric interpretation of
the phase which appears in the scalar product of two coherent vectors.
 In  1994, in Bia\l owie\.{z}a, Askold Perelomov added this
  question on my list \cite{sb6} of other
6 questions referring to coherent states and geometry.
 An explicit answer to this question  for the Riemann sphere
 is given by Perelomov himself (cf. Ref. \cite{per}, p. 63). Earlier,
S.  Pancharatnam \cite{pan,wil} showed that the phase difference between
 the
initial and final state is $<A|A'>=\exp (-i\Omega_{ABC}/2)$, where
$\Omega_{ABC}$ is the solid angle subtended by the geodesic triangle $ABC$
on the Poincar\'e sphere.
 The holonomy of a loop in the projective Hilbert space is twice the
symplectic area of any two-dimensional submanifold whose boundary is the given
loop (see Prop. 5.1 in 
\cite{aa}, where this result is attributed to Aharonov and Anandan).

The main result communicated at this conference is the following: on certain 
manifolds,  the phase $\Phi$ which appears in the scalar product of
two coherent state vectors is twice the symplectic area of the geodesic 
triangle determined by the corresponding points on the manifold and the origin
of the system of coordinates. This result was proved on a restricted class
of manifolds: the compact, homogeneous, simply connected Hodge manifolds, which
are in the same time naturally reductive.  I mention also that
during the Workshop  Martin Bordemann pointed out that the class of manifolds
considered by me consists in
fact only of  the Hermitian symmetric spaces \cite{martinb}. Indeed,
any naturally reductive space with an invariant K\"ahler structure is  locally
 Hermitian
symmetric \cite{van} and simply connectedness implies Hermitian symmetry.
 On the other side, the results of the present paper are still true for  other
manifolds than those considered here. For example, the results are true for
 the Heisenberg-Weyl group \cite{per} as well as for the
noncompact dual of the  complex Grassmann manifold \cite{sb33}.

One  idea of the present paper is to  use  in connection with geodesic
 triangles the generalization  via coherent states  of
the shape invariant
 \cite{blaschke}.  A more precise formulation
on local and global realization 
are given in \cite{sb3}.    Details will be given elsewhere \cite{mars}.

The paper is laid out as follows. In \S 2 the notation on
 holomorphic line bundles and  coherent states is fixed.
\S 3 deals with the shape invariant of Blaschke and Terheggen and its
generalization to coherent states.  The
 results on phases of coherent states are proved in \S 4.
\S 5 presents briefly the calculation of the symplectic area of geodesic
 triangle on  \Gras~ \cite{sb33} using the technique  and notation
from  Ref. \cite{sb2}.

\section{\hspace{-4mm}.\hspace{2mm}
HOLOMORPHIC LINE BUNDLES AND COHERENT STATES}

Let $\tau: L\rightarrow M$ be a holomorphic line bundle over the \Ka~ manifold
 $(M, \omega )$, with the connection $\nabla_L$ compatible  with  the
hermitian metric $h_L$.  With respect  to a holomorphic frame and a
holomorphic coordinate system,
$\nabla_L=\partial + \theta_L+ {\overline{\partial}}$,
$\theta_L=\partial\log {\hat{h}}_L$, and
 $\Theta_L=\overline{\partial}\theta_L$,
where  $\theta_L$ ($\Theta_L$) is the connection (respectively, curvature) 
matrix.

Assume that  $\tau$ is a  prequantum line bundle, i.e.
$\omega =\frac{\sqrt{-1}}{2}\Theta_L= \pi c_1(L)$.  Then $M$ is a  Hodge
 manifold and $L$ is a  positive, ample   line bundle, here taken  very
 ample. The embedding $ \iota :M\hookrightarrow \gcp ^N$
 is assumed  holomorphic and isometric,
 which implies that  it is  \Ka ian: $\omega_M=\iota^*\omega_{FS}$.
Above  $N=\mbox{\rm{dim}}\, \gh -1$, where $ \gh =H^{0}(M,L)=
 \Gamma_{\mbox{\rm{hol}}}(M,L)$. If $ s_i , i=1,...,N+1$ is
 a  basis of global sections, orthonormal with respect with the scalar
 product on $\gh $, then the embedding $\iota $ is given by
\begin{equation}\label{iota1}
\iota (z)=(s_1(z):s_2(z):...:s_{N+1}(z)).
\end{equation}

  Rawnsley's \cite{raw} coherent states are defined as usual:
if  $q\in L\setminus\{ 0\}=L_0,$ is fixed, then
 the evaluation of the section  $ s\in \gh$  determines uniquely the 
coherent vector $e_q\in\gh$ , $s(\tau (q))=
(e_q,s)q$.

 Perelomov's \cite{per} coherent states  are defined by the triplet
 $(G,\pi ,\gh)$, where $G$ is a  Lie group, $\pi$ a  unitary irreducible
representation on the
 complex separable Hilbert space $\gh$.
 Let $e_0\in \gh$ be fixed and $e_g=\pi
(g) e_0$. With the notation
 $\tilde{\psi}=\{e^{i\alpha}\psi|\alpha\in \gr\}, \psi\in\gh$,
$\{e_g\}_{g\in G}$  is a family of coherent vectors, while
$\{\tilde{e}_g\}_{g\in G}$  is a family of coherent states.
If $K=\{k\in G|\pi (k)e_0=e^{i\alpha (k)}e_0\}$, then
$ M=\{\tilde{\pi}(g)\tilde{e}_0|g\in
G\}$ and $M\approx G/K$.   Let
 $\chi$ be a character of $K$. Then, in Perelomov's construction, $L=
 M\times_{\chi}\gc $ is a  $G$-homogeneous line bundle associated
  by the character $\chi$ to the principal $K$-bundle
 $K\rightarrow G\rightarrow M$.
\nopagebreak
 In fact, Perelomov's \cite{per} coherent vectors are
\nopagebreak
\begin{equation}
{\mb e}_{Z,j}=\exp\sum_{{\varphi}\in\Delta^+_n}(Z_{\varphi}F^+_{\varphi})
j ,~~~~\underline{{\mb e}}_{Z,j}=({\mb e}_{Z,j} ,
{\mb e}_{Z,j})^{-1/2}{\mb e}_{Z,j} ,
\label{z}
\end{equation}
 where $\Delta^+_n$ are
 the positive non-compact roots, $ Z\df (Z_ \varphi )
\in {\gc}^D$  are local
coordinates in the maximal neighbourhood
${\cal V}_0 \subset  M $, $F^+_{\varphi} j\neq  0, (F^-_{\varphi} j = 0),
~  \varphi\in\Delta^+_n $, and $j$ is the extremal weight vector
of the representation $\pi$ (see also \cite{sb4}).

\section{\hspace{-4mm}.\hspace{2mm}THE GENERALIZATION TO COHERENT STATES
OF THE SHAPE INVARIANT OF BLASCHKE AND TERHEGGEN}

Let us consider the projection
$\xi: \gc^{n+1}\setminus \{0\}\rightarrow \gcp^{n},~\xi (x)=[x]$.

Let us consider the function: $\Psi_{\cpn} :\gcp^n\times\gcp^n\times
\gcp^n\rightarrow\gc$
\begin{equation}\label{psii}
\Psi_{\cpn} ([x],[y],[z])=\frac{(x,y)(y,z)(z,x)}{||x||^2||y||^2||z||^2},~
x,y,z\in \gc^{n+1}\setminus \{0\},
\end{equation}
where the
scalar product $( x, y)$ in  $\gc^{n+1}$ is linear in the second entry.
 Let us use the notation
$d_C([x],[y])=\arccos \frac{|(x,y)|}{||x|| ||y||}$ for the Cayley distance.

The phase $\Phi$ on a closed loop passing through three states
 in the projective space was considered
by Bargmann \cite{bar}.
Here I correlate this phase   with the shape invariant used by
 Blaschke  and Terheggen  for $\gcp^2$
 \cite{blaschke} and by  Brehm for $\gcp^n$ \cite{brem}. They have
proved that:
$$\Psi_{\cpn} ([x],[y],[z])=\cos a \cos b \cos c\, \exp ({i\Phi_{\cpn}}),$$
where $0\le\Phi_{\cpn} < 2\pi$,
and $a,b,c<\pi /2$  (in order to assure the uniqueness of the geodesic arcs)
are the  sides of the triangle $[x], [y], [z]$:
$a=d_C([y],[z]),~b=d_C([z],[x]),~c=d_C([x],[y]).$

\begin{thm}[Hangan, Masala \cite{hangan}]
{ Given a geodesic triangle with vertices
$[x], [y]$, $ [z]$ in the projective space $\gcp^{n}$, let $S$ be the surface
generated by the geodesic arcs issued from $[x]$ with end-points on the
 geodesic
arc between $[y]$ and $[z]$. Let $\got{I}$ be the integral of the two-form
$\omega$ on $S$. Then
\begin{equation}\label{han}
\Phi_{\cpn} = -2\got{I}+2k\pi, k\in\gz.
\end{equation}}
\end{thm}
\begin{rem}[Hangan, Masala]{ As $\omega$ is closed, we have as a 
consequence of
 Stokes' theorem
that $\got{I}$ does not vary when $S$ is continuously deformed such that its
boundary is fixed.}
\end{rem}
Theorem \ref{sb} below
gives  in the particular case of the complex projective space
 a new proof of the theorem of Hangan and Masala.
 For the complex Grassmann manifold, an explicit  calculation, 
independent of Theorem \ref{sb2} is briefly
presented in \S 5.

Now we shall consider a generalization of the definition (\ref{psii}) for the
line bundle $\tau$ in the context of Rawnsley's coherent states. So, let us
 take $x,y,z\in
M$ and $q,q',q''\in L$ such that $\tau (q)=x, \tau (q')= y, \tau (q'')= z$. A
generalization of (\ref{psii}) is given by the three-point function
$\Psi_{M} :M\times M\times M \rightarrow \gc$  (see also \cite{sb3}) :
\begin{equation}\label{psiii}
\Psi_{M} (x,y,z)=\frac{(e_q,e_{q'})(e_{q'},e_{q"})(e_{q"},e_{q})}
{||e_q||^2||e_{q'}||^2||e_{q"}||^2},
\end{equation}
which is globally defined and does not depend of the representatives in the
fibre.

\begin{thm}\label{sb}
{a). Let $M$ be a compact,  Hodge manifold, admitting the k\"ahlerian 
embedding (\ref{iota1}).
 Then we have the Cauchy
formula:
\begin{equation}\label{cau}
\Psi_{M} (x,y,z)=\Psi_{\cpn} (\iota (x), \iota (y), \iota (z)).
\end{equation}
Let $0\le\Phi_{M}<2\pi$ be the phase
\begin{equation}\label{unghi}
\Psi_{M} (x,y,z)=|\Psi_{M} (x,y,z)|\times e^{i\Phi_{M} (x,y,z)}.
\end{equation}
Then we have
\begin{equation}\label{embpsi}
\Phi_{M} (x,y,z)=\Phi_{\cpn} (\iota (x), \iota (y),
 \iota (z))~ \mbox{\rm{mod}}~
2k\pi ,k\in\gz,
\end{equation}
\begin{equation}\label{cs}
|\Psi_{M} (x,y,z)|=\cos\,a\cos\,b\cos\,c ,
\end{equation}
where
\begin{equation}
a=d_C(\iota (y),\iota (z)), b=d_C(\iota (z),\iota (x)),
 c=d_C(\iota (x),\iota(y)) .
\end{equation}

b). Let us suppose that $M$ is a compact Hermitian symmetric space.
Let us consider that the  points $x,y,z\in M$ are such that  any pair of them
can be joined by a unique geodesic arc, which determine the loop
$\gamma (x,y,z)$ . Then the
 angle $\Phi_{M}$ (\ref{unghi})
  can be given by  an equation of the type
 (\ref{han}), but on the manifold $M$:
\begin{equation}\label{han2}
\Phi_{M} (x,y,z)=-2\int _{\sigma (x,y,z)}\omega_M.
\end{equation}
where  $\sigma (x,y,z)$ is the surface of the geodesic triangle determined
by the points $x,y,z$ or a deformation surface of $\gamma (x,y,z)$}.
\end{thm}
{\it Proof}.
a). Eq. (\ref{cau}) is  proved using the Cauchy formula \cite{sb6}
 for the complex
 two-point functions in  a local representation of sections.
 See also Prop. 4.7 in Ref. \cite{sb3}.

 Let us also consider the real two-point function
\begin{equation}\label{is}
\psi_{M} (x,y)=\frac{|(e_q,e_{q'})|}{||e_q||^{1/2}||e_{q'}||^{1/2}},
 \tau (q)=x, \tau (q')= y.
\end{equation}

The  two-point function verifies locally  the Cauchy relation
\begin{equation}\label{cusi}
\psi_{M} (x,y)= \psi_{\cpn} (\iota (x), \iota (y))
\end{equation}
(see
\S 4.2 in Ref. \cite{sb3} for a more precise formulation).
The functions (\ref{is}) are introduced in eq. (\ref{psiii}), the  Cauchy
relation (\ref{cusi}) is taken into account and eq. (\ref{embpsi}) follows.
Eq. (\ref{cs}) follows if in the Cauchy relation (\ref{cusi}) 
 it  is observed that
\begin{equation}\label{ce}
\psi_{\cpn} (\iota (x), \iota (y))=\cos d_C (\iota (x), \iota (y)).
\end{equation}

b). Now we consider the manifold $M$ to be to be compact Hermitian symmetric.
 So, $M$ is
 a compact, homogeneous, simply connected, naturally reductive, Hodge manifold,
which admits a holomorphic and isometric embedding in a projective
 space.

Let us consider a closed piece-wise smooth curve $\gamma$ in $M$. Because the
 manifold
is     Hodge and simply connected, we are  under the
conditions of
 Thm. 2.2.1 in \cite{kostant}. From  \cite{kostant}
 we need only the expression (1.8.3) of the 
parallel transport function. However,
 in order to put in accord the notation from \cite{kostant}
 with our notation, we repeat some parts of the proof. The same notation is
 used also in Theorem \ref{gr} below.
  Let $P_{\gamma}:L_p\rightarrow L_p$ the parallel transport along $\gamma$.
 Then $P_{\gamma}(s)=Q_{\gamma}s$.
Let the notation $A_L=i\theta_L$. The  
 scalar parallel transport function 
$Q(\gamma )=\exp i\beta$ is calculated with the Stokes' formula,
where the phase $\beta$ is
\begin{equation}\label{ber}
\beta =\oint_{\gamma}A_L=\oint_{\gamma}i\theta_L=  \int_{\sigma}dA_L.
\end{equation}
Here $\sigma$ is a surface of deformation of $\gamma$. We recall that
 $\gamma :I\rightarrow M$ is homotopic to a point if there is a rectangle
 $R=[a,b]\times [c,d]$ in the plane and a piece-wise smooth parametrization
 $\rho: I\rightarrow \stackrel{\circ}{R}$ of the boundary
 $\stackrel{\circ}{R}$ of $R$ oriented counter-clockwise such that
 $\sigma\circ\rho =
 \gamma$. Such a map defines an oriented surface with $\stackrel{\circ}{R}$
 oriented counter-clockwise as boundary and is called surface of deformation of
 $\gamma$.

 Considering the  (positive) line bundle $L$ over $M$,  we have (also cf.
eq. (2.2.1) in \cite{kostant}) \begin{equation}\label{bt}
\beta =i\int_{\sigma}d\theta_L=i\int_{\sigma} \Theta_{L}
=2\int_{\sigma}\omega_M,
\end{equation}
\begin{equation}\label{scalar}
Q(\gamma )=\exp ({i\beta})=\exp({2i\int_{\sigma}\omega_M }).
\end{equation}

Now  a closed path $\gamma :x\rightarrow y\rightarrow z\rightarrow x$ in
$M$ is considered.
Because the manifold $M$ is naturally reductive, the coherent states realise
parallel
transport on geodesics (cf. Remark 3 in the second Ref. \cite{sb4} and
 Remark 1 in the third Ref. \cite{sb6}). Taking as auto-parallel 
section $s$  along the piece-wise smooth curve $\gamma$
in the formula of the parallel transport
  $P_{\gamma}(s)=Q_{\gamma}s$  the normalized coherent state vector 
$||e_q||^{-1}e_q$, the  holonomy $\beta = \Phi (x,z,y)$ on the
  geodesic path  $\gamma =  \gamma (x,y,z,x)$ is
$\beta = -\Phi (x,y,z).$
So, eq. (\ref{han2}) is proved. More details will be given elsewhere
\cite{mars}. \hfill$\pre$

\section{\hspace{-4mm}.\hspace{2mm}PHASES AND COHERENT STATES}
\begin{thm}\label{sb2}{Let $M$ be a compact Hermitian symmetric manifold.
 Let $(L,h_L,\nabla_L)$ be a
homogeneous line bundle supposed to be very ample. Let us consider on the
manifold of coherent states $M$ the Perelomov's coherent vectors (\ref{z}) in a
local chart, corresponding to the fundamental representation $\pi$.
Let us consider the points $Z, Z'\in {\cal V}_0\subset M $ such that
$0, Z, Z'$ is a geodesic triangle.
 Then the phase $\Phi_M$ defined by the relation
\begin{equation}\label{pufi}
(\underline{{\mb e}}_{Z'},\underline{{\mb e}}_{Z})=
|(\underline{{\mb e}}_{Z'},\underline{{\mb e}}_{Z})|\exp({i\Phi_M (Z',Z)})
\end{equation}
is given by twice  the integral of the  symplectic two-form on the
surface  $\sigma (0,z,z')$ of the geodesic triangle
 $\gamma (0,Z,Z')\subset M$
 \begin{equation}\label{pha}
 \Phi_M (Z',Z)=2\int_{\sigma (0,Z,Z')}\omega_{M}.
 \end{equation}

Also
\begin{equation}\label{dcc}
|(\underline{{\mb e}}_{Z'},\underline{{\mb e}}_{Z})|=
|(\underline{{\mb e}}_{\iota (Z')},\underline{{\mb e}}_{\iota (Z)})|=
\cos d_C(\iota (Z'), \iota (Z)).
\end{equation}}
\end{thm}
\begin{rem}{If the line bundle is not very ample, then an integer $m$ appears
in
front of the integral in the eq. (\ref{pha}), corresponding to the power
$L^m$ for which the ample line bundle $L$ becomes very ample.}
\end{rem}
{\it  Proof of the Theorem \ref{sb2}.\/} The theorem is a direct
 consequence of
Theorem \ref{sb} and of the remark
that the Berry phase \cite{wil} is the opposite of the Bargmann phase 
\cite{bar}.
Indeed, let us take the points $0, Z, Z'$ in Theorem \ref{sb2} to
correspond to the
points $x, y, $ respectively $z$ in Theorem \ref{sb}. Then the three-point
function $\Psi _M(0,Z,Z')$ (\ref{psiii}) becomes  the complex-valued
 two-point function
$\frac{(\underline{{\mb e}}_{Z},\underline{{\mb e}}_{Z'})}
{||\underline{{\mb e}}_{Z}||^{1/2}||\underline{{\mb e}}_{Z'}||^{1/2}}$
and $\Phi_M (0,Z',Z)=-\Phi_M (0,Z,Z')$ in eq. (\ref{unghi})
 is denoted simply $\Phi_M (Z,Z')$ in eq. (\ref{pufi}).
Eq. (\ref{dcc}) is nothing else than eq. (\ref{ce}).\hfill \pre

\section{\hspace{-4mm}.\hspace{2mm}{ILLUSTRATION ON COMPLEX GRASSMANN MANIFOLD
\Gras}}

\begin{thm}\label{gr}{Let $z,z'\in{\cal V}_0\subset\Gras$ be described by the
Pontrjagin coordinates $Z, Z'$. Let $\gamma (0,z,z')$ be the geodesic
triangle obtained by
joining $0, z, z'$. Then the symplectic area of the surface $\sigma (0,z,z')$
 of the geodesic triangle $\gamma (0,z,z')$ is given by
 \begin{equation}\label{ar}
 {\got I}=\int_{\sigma (0,z,z')}\omega
 =\frac{1}{4i}\log\,\frac{\det(1+ZZ'^+)}
 {\det(1+Z'Z^+)}.
 \end{equation}}
 \end{thm}
{\it Proof\/} We apply the Stokes' formula (\ref{ber}), take into
account  eqs. (\ref{bt}), (\ref{scalar}) and the relation
$dA_L=2\omega .$
 Here $L$ is the dual of the tautological (universal) line
bundle on the Grassmann manifold.
The connection is
$A_L=i \mbox{\rm{Tr}}\,[dZ\,Z^+(1+ZZ^+)^{-1}]$  \cite{sb4}.
In the relation $A_L=i\theta_L$,  the connection matrix $\theta_L$
corresponds to the hermitian
metric on the dual of the  tautological line bundle on the Grassmann manifold
$\hat h_{L}(Z)=\det (1+ZZ^+)^{-1}$.
The Berry connection \cite{wil} which corresponds to $A_L$ is
\begin{equation}
A_B=\frac{i}{2} \mbox{\rm{Tr}}[(dZ\,Z^+-Z\,dZ^+)(1+ZZ^+)^{-1}].
\end{equation}
The  corresponding two-form  on \Gras~is
\begin{equation}\label{ogras}
\omega=\frac{i}{2}\mbox{\rm{Tr}}[dZ(1+Z^+Z)^{-1}\wedge dZ^+(1+ZZ^+)^{-1}].
\end{equation}

The calculation is  long. I indicate here only the main steps. Details will be
given elsewhere \cite{sb33}.

a). Firstly, let $z\in {\cal V}_0$. The explicit expression of the geodesic
starting at $0\in {\cal V}_0\subset \Gras$  with $\stackrel{.}{Z}(0)=B$
is
$Z(t)=B\tan \frac{\sqrt {B^+B}t}{\sqrt {B^+B}}$, where $\dot{Z}(0)=B$ (cf.
\cite{sb2}).

The conditions that the points $0, z, z_0\in {\cal V}_0$ to belong to  the
 same geodesics are \cite{sb33}
\begin{equation}\label{conditia}
Z_0Z^+=ZZ_0^+; Z_0^+Z=Z^+Z_0.
\end{equation}
b). The integral on $\gamma (0,z_1,z_2,0)$ is calculated firstly on the
geodesic
arc joining $z_1,z_2$.  The situation is reduced to that on calculating the
 integral
on the geodesic joining the points $0,z$. A linear
 fractional transformation which sends $z_1\rightarrow 0$ has
the expression \cite{sb2}:
$Z'(Z)=(AZ+B)(CZ+D)^{-1}$
where
$$A\!=\!(\oo + Z_1Z_1^+)^{-1/2},\!
                        B\!=\!-(\oo + Z_1Z_1^+)^{-1/2}Z_1,\!
                        C\!=\! (\oo + Z_1^+Z_1)^{-1/2}Z^+_1,\!
                        D\!=\!(\oo + Z_1^+Z_1)^{-1/2}$$
\begin{equation}
Z'(Z)  =  (\oo + Z_1Z^+_1)^{-1/2}(Z-Z_1)(\oo + Z^+_1Z)^{-1}
(\oo + Z^+_1Z_1)^{1/2}.\label{zz}
\end{equation}

When $Z_1\rightarrow 0$, the point $Z_2$ becomes $Z_I=Z'(Z_2)$.
So, we have to calculate
\begin{equation}\label{cal}
I=i\int_{Z_1}^{Z_2}\mbox{\rm{Tr}}\,[dZZ^+(1+ZZ^+)^{-1}]
\end{equation}
and make the quasi-linear change of variables. We need the formulas:
\begin{equation}\label{lll}
\begin{array}{l}
dZ  = (A-Z'C)^{-1}dZ'(B^+Z'+D^+)^{-1},
Z^+  =  (D^+Z'^+-B^+)(A^+-C^+Z'^+)^{-1},\\
1+ZZ^+  = (-Z'C+A)^{-1}(1+Z'Z'^+)(C^+Z'^+-A^+)^{-1}.
\end{array}
\end{equation}

Eq. (\ref{cal}) becomes
\begin{equation}\label{i12}
I=i\int_0^{Z_I}\mbox{\rm{Tr}}\,[dZ'(1-Z_1^+Z')^{-1}(Z'^++Z_1^+)(1+Z'Z'^+)^{-1}].
\end{equation}

c). The condition (\ref{conditia}) under the fractional transformation
 becomes
 \begin{equation}
Z_IZ'^+= Z'Z_I^+; Z'^+Z_I=Z^+_IZ'.
\end{equation}  \label{ecr}
The last equation has the solution
$Z'=\frac{|Z_{11}|}{|Z_{I11}|}Z_I$ \cite{sb33}.
Introducing the last expression in eq. (\ref{i12}) it is obtained
\begin{equation}\label{vul}
I=i\int_0^{r_0}dr\mbox{\rm{Tr}}\,[(1-Ar)^{-1}(Br+A)(1+Br^2)^{-1}],
\end{equation}
where
$A=\frac{Z^+_1Z_I}{|Z_{I11}|}, |Z_{1,11}|=r,B= \frac{Z_I^+Z_I}{|Z_{I11}|^2},
|Z_{I11}|=r_0.$
But
$$(1-Ar)^{-1}(Br+A)(1+Br^2)^{-1}= Br(1+Br^2)^{-1}+(1-Ar)^{-1}A.$$
With the formula
$\frac{d}{dx}\log\det U=\mbox{\rm{Tr}}(U^{-1}\frac{\partial U}{\partial x})$,
the integral (\ref{vul}) becomes successively
\begin{equation}I=i\log\det \frac{(1+Br^2_0)^{1/2}}{1-Ar_0}=
\frac{i}{2}\log\det\frac{(1+Z_I^+Z_I)}{(1-Z_1^+Z_I)^2}.
\end{equation}
But
\begin{eqnarray}
\nonumber 1-Z_1Z^+_I &\! =\! & (1+Z_1Z_1^+)^{1/2}
 (1+Z_1Z_2^+)^{-1}(1+Z_1Z_1^+1)^{1/2},\\
\nonumber (1+Z_IZ_I^+)^{-1} & = &
\nonumber (1+Z_1Z_1^+)^{-1/2}(1+Z_1Z_2^+)^{-1/2}\times\\
\nonumber & & (1+Z_2Z_2^+)^{-1}(1+Z_2Z_1^+)(1+Z_1Z_1^+)^{-1/2},
\end{eqnarray}
\begin{equation}
I=\frac{i}{2}\log\frac{\det (1+Z_2Z_2^+)\det (1+Z_2Z_1^+)}
{\det (1+Z_1Z_1^+)\det (1+Z_1Z_2^+)}.
\end{equation}
Taking the particular values $0,Z$ and $Z'$
 for $Z_1,Z_2$ in the last expression, eq.
(\ref{ar}) is proved because $I=2{\got I}$.\hfill\pre

\begin{rem}{Equation (\ref{ar}) contains as particular case the projective
 space
and the sphere. The expression for the sphere can be found in \cite{per}.
 Note that
in the conventions of this talk, the two-form on the sphere is
$\omega =\frac{i}{2}\frac{dz\wedge d\overline{z}}{(1+|z|^2)^2}$.
This gives for the sphere of radius 1 the area $\pi$. This also explains the
difference with Panchartnam's formula.}
\end{rem}

{\bf Acknowledgments}

The author expresses his thanks to the organizers of the XVII-th Workshop on
 Geometrical
Methods in Physics 
 for the inviting him at  the Workshop. 
Suggestions from   Professors 
  M. Bordemann, P. Delano\"e, C. Gheorghe, T. Hangan, L. Ornea, B. Rosenfel'd,
 G. Tuynman, M. Schlichenmaier
and L. Vanhecke are kindly acknowledged.

\end{document}